\documentclass[final]{siamart171218}

\usepackage{soul}
\usepackage{arydshln}
\usepackage{mathtools}
\usepackage{lipsum}
\usepackage{amsfonts}
\usepackage{graphicx}
\usepackage{epstopdf}
\usepackage{algorithmic}
\ifpdf
  \DeclareGraphicsExtensions{.eps,.pdf,.png,.jpg}
\else
  \DeclareGraphicsExtensions{.eps}
\fi

\newsiamremark{remark}{Remark}
\newsiamremark{hypothesis}{Hypothesis}
\crefname{hypothesis}{Hypothesis}{Hypotheses}
\newsiamthm{claim}{Claim}
\newsiamthm{example}{Example}
\newsiamthm{counterexample}{Counterexample}
\headers{Local Controllability of Quadratic Systems}{Moise R. Mouyebe, Anthony M. Bloch}

\title{On the Local Controllability of a Class of Quadratic Systems\thanks{Submitted to the editors.}}
\author{Moise R. Mouyebe\thanks{Department of Mathematics, University of Michigan, Ann Arbor, MI 48109, USA
  \email{\{mmouyebe,~abloch\}@umich.edu}.}
\and Anthony M. Bloch\footnotemark[2]}

\usepackage{amsopn}

\makeatletter
\newcommand*{\addFileDependency}[1]{% argument=file name and extension
  \typeout{(#1)}% latexmk will find this if $recorder=0 (however, in that case, it will ignore #1 if it is a .aux or .pdf file etc and it exists! if it doesn't exist, it will appear in the list of dependents regardless)
  \@addtofilelist{#1}% if you want it to appear in \listfiles, not really necessary and latexmk doesn't use this
  \IfFileExists{#1}{}{\typeout{No file #1.}}% latexmk will find this message if #1 doesn't exist (yet)
}
\makeatother

\newcommand*{\myexternaldocument}[1]{%
    \externaldocument{#1}%
    \addFileDependency{#1.tex}%
    \addFileDependency{#1.aux}%
}
\ifpdf
\hypersetup{
  pdftitle={An Example Article},
  pdfauthor={D. Doe, P. T. Frank, and J. E. Smith}
}
\fi

\myexternaldocument{ex_supplement}

\begin{document}

\maketitle

% REQUIRED
\begin{abstract}
  The local controllability of a rich class of affine nonlinear control systems with nonhomogeneous quadratic drift and constant control vector fields is analyzed. The interest in this particular class of systems stems from the ubiquity in science and engineering of some of its notable representatives, namely the Sprott system, the Lorenz system and the rigid body among others. A necessary and sufficient condition for strong accessibility reminiscent of the Kalman rank condition is derived, and it generalizes Crouch's condition for the rigid body. This condition is in general not sufficient to infer small-time local controllability. However, under some additional mild assumptions  local controllability is established. In particular for the Sprott and Lorenz systems, sharp conditions for small-time local controllability are obtained in the single-input case.
\end{abstract}

% REQUIRED

% REQUIRED

\section{Introduction}

The use of linear algebraic methods is  pervasive in science and engineering. This can be attributed to its very intuitive and widely accessible language that is also well suited for computer implementation. In the field of control systems in particular, these methods have been remarkably successful in the analysis and design of linear systems. As a result, even when dealing with nonlinear systems, scientists and control engineers alike are more likely to start their analysis by considering as a first step the linearization of the system. However, a linearization even though convenient, is only an approximation that does not capture the full complexity of the nonlinearity embedded in the dynamics of most real life systems. Hence, roughly since the seventies, pioneering work by Brockett \cite{brockett1972algebraic}\cite{brockett1973lie}\cite{brockett1973lie-B}, Haynes and Hermes \cite{haynes1970nonlinear}, Hermann \cite{hermann1963accessibility}, Hermes  \cite{hermes1974local}\cite{hermes1976local}\cite{hermes1978lie}, Lobry \cite{lobry1974controllability}, Sussmann and Jurdjevic \cite{sussmann1972controllability}, Krener \cite{krener1974generalization}, Hermes and Kawski \cite{hermes1986local} to name a few has been instrumental  in anchoring the foundation of the field of nonlinear control. A key development in the field was  the general theorem on local controllability by Sussmann \cite{sussmann_Gen_th_on_Local_ctrb} in the late eighties followed in the early nineties by an improvement due to Bianchini and Stefani \cite{bianchini1993controllability}, grachev
and Gamkrelidze \cite{agrachev1993local} among others, all of whom give general sufficient conditions for local controllability. There is still to 
 our knowledge no sharp condition, namely necessary and sufficient,   that is  applicable to generic nonlinear systems. This is certainly not because of the lack of interest in research in the field, but it's rather due to its inherent complexity. The  tools of differential geometry and  Lie algebraic theory as in the work of \cite{brockett2014early}
 are still of interest today, and in this work we are able to use these ideas and the theory developed in 
 \cite{crouch1984spacecraft} to extend the theory and apply to some new systems of interest. 
 
 The first question one must answer in the analysis of any control system is that of whether the system is controllable. For linear systems, the answer is given by the Kalman rank condition. In nonlinear setting however, the answer is not as straightforward as computing the rank of some subspace, and sometimes it can be quite hard \cite{sontag1988some}\cite{kawski1990complexity} that it's worth considering a weaker notion than controllability such as accessibility which in fact is necessary for any form of controllability. In this paper, we analyze  the local controllability of a specific class of quadratic affine control systems. The interest in this class of systems  comes from the ubiquity in science and engineering of some of its notable members among which the Sprott system \cite{sprott2010elegant}\cite{kvalheim2021families} which is a tractable version of the repressilator \cite{elowitz2000synthetic}\cite{bucse2009existence}\cite{buse2010dynamical}, a model of synthetic genetic regulatory network consisting of a ring oscillator and widely used in biology; the rigid body \cite{bloch2015nonholonomic}\cite{arnol2013mathematical}\cite{marsden1997introduction}  which is fundamental in mechanical and aerospace engineering \cite{crouch1984spacecraft}; and the Lorenz system \cite{lorenz1963deterministic} which models dissipative fluid flow with vast application in science and engineering \cite{haken1975analogy}\cite{cuomo1993circuit}\cite{poland1993cooperative}. The uncontrolled dynamics of these three systems in particular reveal quite disparate behaviors. For instance, the Sprott system hints at chaotic behavior in the non-dissipative regime \cite{sprott2010elegant}, the rigid body is conservative and Poisson-Hamiltonian \cite{marsden2013introduction}, and the Lorenz system is dissipative and chaotic \cite{lorenz1963deterministic}. Meanwhile, their controlled dynamics have remarkable similarities since they can be generically analyzed within the framework of the class of quadratic affine systems we're considering here. Prior work along these lines came  among others from Jurdjevic and Kupka \cite{jurdjevic1985polynomial}, Brunovsky \cite{brunovsky1974local} who studied the local controllability of  polynomial control systems with an emphasis on those with an odd degree of homogeneity for which small-time local controllability is equivalent to the Lie Algebra Rank Condition. Melody et al \cite{melody2003nonlinear}  studied the case of homogeneous quadratic polynomial drift and emphasized the use of congruence transformation as a mean to neutralize bad brackets. More recently, Aguilar \cite{aguilar2012local} studied the local controllability of a restrictive class of homogeneous quadratic affine control systems in which he remarkably exploited the graph theoretic calculation tool of labeled trees introduced by Butcher \cite{butcher2016numerical} to obtain a sufficient condition for small-time local controllability in terms of higher order variations. His work was subsequently generalized by Krastanov \cite{krastanov2023small}.
 
 In this work, we leverage the fact that the free dynamics of the systems under consideration is afforded by non-necessarily homogeneous quadratic  vector fields to aggregate them into a larger class of quadratic affine  control systems. This allows us to give a uniform treatment of their shared local controllability properties. First, we derive a necessary and sufficient condition for strong accessibility for this larger class of quadratic systems emulating the Kalman rank condition, thus making it suitable for computer implementation as a first check for controllability, and we show that our condition recovers Crouch's \cite{crouch1984spacecraft} for the rigid body. Next, we address the small-time local controllability guarantee of our condition. Finally, we give explicit necessary and sufficient conditions for single-input small-time local controllability for the Sprott and Lorenz systems.

\subsection{Setup and Notation}
Throughout this paper, we study  affine control systems in $\mathbb{R}^n$ of the form
\begin{align}\label{class_S_Ctrl}
  \Sigma_k: \quad\quad  \dot{x} = f_0(x)+u_1f_1+\cdots +u_{n-k}f_{n-k}, \quad x(0)=x_0; \quad 1\le k \le n-1
\end{align}
with control input $u = (u_1,\cdots,u_{n-k})$ constrained in a bounded set $U$ satisfying Aff$(U) = \mathbb{R}^{n-k}$, where Aff$(U)$ denotes the affine hull of $U$, namely the set of all finite linear combinations $\sum{\alpha_i\nu_i}$ with $\nu_i \in U$, $\alpha_i \in \mathbb{R}$, and $\sum{\alpha_i}=1$. Admissible controls are Lebesgue integrable functions $u(\cdot): [0,T]\rightarrow U$ for some positive real $T$. The control vector fields $f_i,~ 1\le i\le n-k$ are assumed to be linearly independent and constant throughout. The drift $f_0$ is a quadratic vector field in the class 
\begin{align}\label{class_S}
   \mathcal{C}: \quad Lx + \sum_{\nu=1}^{n}{Q_{\nu}(x_{\nu+1[n] },x_{\nu+2[n] })}\frac{\partial}{\partial x_{\nu}} 
\end{align}
where $L$ is a linear operator in $\mathbb{R}^n$, and  the $Q_{\nu}$'s are homogeneous  quadratic polynomials of the form $Q_{\nu}(x_i,x_j) = a_{\nu}x_i^2 + b_{\nu}x_j^2+c_{\nu}x_ix_j$, for $1\le {\nu}\le n$. \emph{$[n]$ stands for the congruence modulo $n$.}

Let $a=(a_1,\cdots,a_n)^T$, $b=(b_1,\cdots,b_n)^T$, $c=(c_1,\cdots,c_n)^T$, and consider the two permutation matrices $P_2$ and $P_3$ which act on the vector $x\in \mathbb{R}^n$ by  cyclically permuting its elements to move $x_2$ and $x_3$ to the first place respectively. More specifically,  $P_2 x = (x_2,x_3,\cdots,x_n,x_1)^T ,$ and $P_3 x = (x_3,x_4,\cdots,x_1,x_2)^T$. Let $v \in \mathbb{R}^n$, we denote by  $\Delta_v$ its diagonal matrix. If $v$ has nonzero components, $\Delta_{\frac{1}{v}}$ will denote the inverse of $\Delta_v$. Finally, for $u,v \in \mathbb{R}^n$, $u \odot v$ will denote their Hadamard product. We define the following quantities:

\begin{align}\label{Psi-Phi-definition}
\begin{split}
    \Psi_{a,b,c}(u,v) &=2a\odot P_2 u \odot P_2 v+2b\odot P_3 u \odot P_3 v + c \odot P_2 u \odot P_3 v + c\odot P_2 v \odot P_3 u\\
    \Phi_{a,b,c}(u) &= \frac{1}{2}\Psi_{a,b,c}(u,u)
\end{split}
\end{align}
With the notation above the drift,  in (\ref{class_S_Ctrl}) takes the form $f_0(x) = Lx + \Phi_{a,b,c}(x)$ where the origin is an equilibrium point, namely $f_0(0)=0$.

Given an admissible control $u(\cdot)\colon [0,T]\rightarrow U$, an admissible trajectory for $u(\cdot)$ is an absolutely continuous curve $x(\cdot)\colon [0,T]\rightarrow \mathbb{R}^n$ that satisfies (\ref{class_S_Ctrl}) for almost all $t\in [0,T]$. Given $p,q\in \mathbb{R}^n$, we say that $q$ is reachable from $p$ in time $T$ if there exists an admissible trajectory such that $x(0)=p$ and $x(T)=q$, and we  denote by $\mathcal{R}(T,p)$ the set of such points. Define $\mathcal{R}(\le T,p)\coloneqq \underset{0\le t\le T}{\cup}\mathcal{R}(t,p)$.

$\Sigma_k$ is said to be \emph{accessible} or has the \emph{accessibility property} (AP) from $p$ if  the interior of $\mathcal{R}(\le T,p)$ is  nonempty; and \emph{locally strongly accessible} from $p$ if the interior of $\mathcal{R}(T,p)$ is nonempty for every $T>0$.   $\Sigma_k$ is said to be \emph{small-time locally controllable} (STLC) from $p$ if $p$ is in the interior of $\mathcal{R}(\le T,p)$ for all $T>0$. Obviously small-time-local controllability and local strong accessibility   imply accessibility. 

The \emph{accessibility algebra} of the affine control  system $\Sigma_k$ is defined  as the smallest Lie algebra of vector fields on $\mathbb{R}^n$ that contains $\mathbf{f}=\{f_0,f_1,\cdots,f_{n-k}\}$, whereas  the \emph{strong  accessibility algebra} $\mathcal{C}_0$ is defined 
as the smallest subalgebra containing $\hat{\mathbf{f}} = \{f_1,\cdots,f_{n-k}\}$ and such that $[f_0,X] \in \mathcal{C}_0$ for all $X \in \mathcal{C}_0$. The \emph{strong accessibility distribution} $C_0$
of $\Sigma_k$   is defined as the distribution generated by the vector ﬁelds in $\mathcal{C}_0$, namely $C_0(p) =  \big\langle X(p) ~\rvert X \in \mathcal{C}_0 \big\rangle$. Lastly, the family $\mathbf{f}$ is said to satisfy the \emph{Lie Algebra Rank Condition}(LARC) at $p$ if the distribution  $L(\mathbf{f})(p)$ associated to the Lie algebra of vector fields generated by the elements of $\mathbf{f}$ spans the whole tangent space. More details can be found in \cite{bloch2015nonholonomic}\cite{nijmeijer1990nonlinear}\cite{sussmann_Gen_th_on_Local_ctrb} for instance and the references therein. 

The Lie bracket $[f,g]$  of the vector fields $f$ and $g$  is given by the formula $[f,g](x) = Dg(x)f(x)- Df(x)g(x)$. We also use the notation $ad_{f}g= [f,g]$, and  recursively define $ad_{f}^{k+1}g = [f,ad_f^kg]$ for $k\ge 0$ with the convention that $ad_{f}^{0}g = g.$ Furthermore, we'll  use  $Span\{ u_1,\cdots, u_s\}$  to denote the subspace spanned by the family $\{u_1,\cdots, u_s\}$.  The sets $\{\mathbf{e}_1,\cdots,\mathbf{e}_n\}$ and $\{\mathbf{E}_{ij}\}_{1\le i,j \le n}$ will respectively denote the canonical bases of $\mathbb{R}^n$ and that of the space Mat$_n(\mathbb{R})$ of square matrices with coefficients in the field $\mathbb{R}$.  Finally,  constant vector fields will  sometimes be safely regarded as ordinary vectors in $\mathbb{R}^n$ in an effort to keep the notation to a minimum.

\subsection{Notable Examples of Systems in this Class}\label{subsection:some_Examples}

\paragraph{The Sprott System \emph{\cite{sprott2010elegant}\cite{kvalheim2021families}}}
This system was originally considered in \cite{sprott2010elegant} in the undamped ($\mu = 0$) regime as an example of "elegant" quadratic chaotic system. The second author of this work will later consider the generalization by adding dissipation ($\mu \ne 0$) to the original system. See \cite{kvalheim2021families} for a rigourous treatment of the existence of period orbits and global bifurcations in this system. It's also worth mentioning that this system can be considered as a tractable version of the repressilator \cite{elowitz2000synthetic} which is  a  model of a synthetic genetic regulatory network consisting of a ring oscillator of odd degree with wide applications in biology.

$L = -(\mu I_3+P_3) = \begin{bmatrix}
    -\mu&0&-1\\-1&-\mu&0\\0&-1&-\mu
\end{bmatrix}$, \quad $\mu \in \mathbb{R}$,\quad  $b  =  \mathbf{0}=(0,0,0)^T =c~$, and $a = \mathbf{1}=(1,1,1)^T$.
The Sprott system dynamics reads:
\begin{align}\label{eq:SprottSystem-Delta-dynamics}
    \dot{x} = Lx + a\odot P_2x \odot P_2x= Lx +\Phi_{\mathbf{1},\mathbf{0},\mathbf{0}}(x)
\end{align}

\paragraph{The Rigid Body \emph{\cite{marsden2013introduction}}}\label{par:RB}
$L = 0$, \quad $a = \mathbf{0} =b$. Let $\xi=(\xi_1,\xi_2,\xi_3)^T \in \mathbb{R}^3$ with $\xi_i>0$ for $1\le i\le 3$ such that the following relation $\Delta_c = \Delta_{\frac{1}{\xi} }(P_2-P_3)(\xi)$ holds for some $c$, then $c\odot P_2x \odot P_3x = \Delta_{\frac{1}{\xi}}S(x)\Delta_{\xi}x$,  where $S(x) = \hat{x} = \begin{bmatrix}
    0&x_3&-x_2\\-x_3&0&x_1\\x_2&-x_1&0
\end{bmatrix}$.  We recover the  rigid body dynamic with inertia matrix $\Delta_{\xi}$ expressed in a suitable basis, namely the column space of the principal inertia matrix.
\begin{align}\label{eq:RB-Delta-dynamics}
        \dot{x} =\Delta_{\frac{1}{\xi}}S(x)\Delta_{\xi}x = \Phi_{\mathbf{0},\mathbf{0},c}(x)
\end{align}

\paragraph{The Lorenz System \emph{\cite{lorenz1963deterministic}}}
$L = \begin{bmatrix}
    -\sigma & \sigma & 0\\
    \rho & -1& 0\\
    0 &0&-\beta
\end{bmatrix}$, 
\quad $\beta,\sigma,\rho >0~$,  
$a = \mathbf{0} = b~$, and $c = (0,-1,1)^T$, which yields the dynamics
\begin{align}\label{eq:Lorenz63-Delta-dynamics}
        \dot{x} = Lx + c\odot P_2x \odot P_3x = Lx + \Phi_{\mathbf{0},\mathbf{0},c}(x)
\end{align}
From now on, we'll drop the explicit dependence of $\Psi$ and $\Phi$ on the vectors $a,b$ and $c$. More specifically, we'll write $\Psi$ for $\Psi_{a,b,c}$ and $\Phi$ for $\Phi_{a,b,c}$ in an effort to further ease the notation.

\section{Necessary and Sufficient Condition for Strong Accessibility}\label{NS-cond-for-strong-local-accessibility}

\begin{figure}[H]
    \centering
    \includegraphics[width=.40\textwidth]{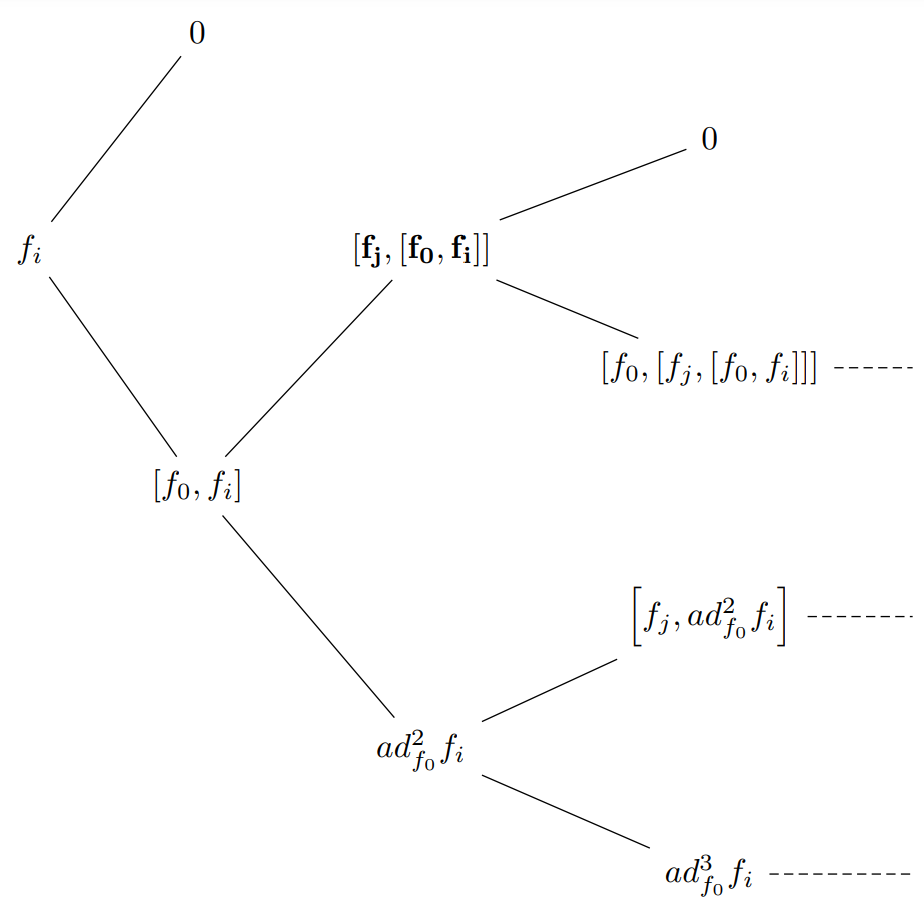} % Image using textwidth to define size
    \caption{Snippet of the strong accessibility algebra generators of the control system $\Sigma_k$. The drift vector field $f_0$ belongs to the class  $\mathcal{C}$ in (\ref{class_S}) and the control vector fields $f_i,f_j,f_l,\cdots$ are constant, $1\le i,j,l,\cdots \le n-k$. The boldface brackets are nonzero constant brackets, namely they don't depend on the state. The next such bracket (not represented in this snippet) is the order $4$ bracket $\mathbf{[f_l,[f_j,ad_{f_0}^2f_i]]}$ that branches out and upward from the bracket $[f_j,ad_{f_0}^2f_i]$. See Table  \ref{long-Brackets} for the values of these brackets at the origin. The values at the origin for the remaining brackets in this figure are as follows: $\mathbf{[f_j,[f_0,f_i]]}=-\Psi(f_i,f_j),~[f_0,[f_j,[f_0,f_i]]](0)=L\Psi(f_i,f_j),~$ and more generally we have $ad_{f_0}^lf_i(0)=(-1)^lL^lf_i,~$ for $l\in \mathbb{N}\cup \{0\}$. } 
    \label{Bracket-Tree}
 \end{figure}

We start by recalling a few important propositions whose proofs can be found in  \cite{nijmeijer1990nonlinear} for instance. 
The first one gives a precise characterization of the strong accessibility algebra $\mathcal{C}_0$ in terms of iterated brackets, and Figure \ref{Bracket-Tree} shows a pictorial representation of this very important result  for the system $\Sigma_k$. Meanwhile, the second one gives a sufficient condition for local strong accessibility.

\begin{proposition}\label{Accessibility-Algebra}
    Every element of $\mathcal{C}_0$ is a linear combination of the iterated Lie brackets of the form $ad_{X_l}\cdots ad_{X_1}f_j$, where $X_i \in \big\{f_0,f_1,\cdots f_{n-k}\big\}$, $1\le j \le n-k$, and $l \ge 0$ with the convention that the case $l=0$ yields $f_j.$
\end{proposition}

\begin{proposition}[Theorem 3.21 in \cite{nijmeijer1990nonlinear}]\label{LocStrgAccSuffCond}
    Consider system (\ref{class_S_Ctrl}), and suppose that $\dim C_0(x_0) = n$, then the system is locally strongly accessible from $x_0$.
\end{proposition}

We have the following brackets:
\begin{align}\label{Thebrackets}
        [f_i,f_0](p) = Lf_i + D\Phi(p)f_i,~ [[f_0,f_i],f_j] = \Psi(f_i,f_j),~ \mbox{and~} [[f_0,f_i],f_i]= 2\Phi(f_i)
\end{align}
The following proposition relates the symmetric bilinear map $\Psi$ defined in (\ref{Psi-Phi-definition}) to the second derivative of the drift $f_0$, then it gives a useful polarization identity that connects the maps $\Phi$ and $\Psi$.
\begin{proposition}\label{polarization-prop}
    Let $f_0$ be the drift in $\Sigma_k$. The following properties hold for any pair of vectors $u,v \in \mathbb{R}^n$:
    \begin{enumerate}
        \item $D^2f_0(u,v) = \Psi(u,v)$
        \item $\Psi(u,v) = \Phi(u+v)-\Phi(u)-\Phi(v)$
    \end{enumerate}
\end{proposition}
\begin{proof}
    The first property comes from a direct computation of the bracket $[[f_0,u],v]$ and using (\ref{Thebrackets}). The polarization identity comes from the symmetric bilinearity of the map $\Psi$ in conjunction with the relation $\Psi(w,w)=2\Phi(w)$ for any $w$ from (\ref{Psi-Phi-definition}).
\end{proof}

\paragraph{The Accessibility Distribution}\label{FundamentalSubspaces}
We inductively define the following fundamental subspace $S_{\lambda},~ \lambda \in \mathbb{N}\cup \{0\}$ as follows:
\begin{align}\label{Subspace:S_lambda}
\begin{split}
  S_0 \coloneqq Span\{ f_1,\cdots,f_{n-k}\}, & \quad\quad S_{\lambda+1} \coloneqq S_{\lambda} + Span\{ L\omega, \Phi(\omega)~\rvert~ \omega \in S_{\lambda}\}
\end{split}
\end{align}

The next couple lemmas are essential in characterizing the accessibility distribution of the control system $\Sigma_k$.

\begin{lemma}\label{C_0-Lemma1}
    $S_{\lambda} \subseteq C_0(0)$ for all $\lambda \ge 0.$
\end{lemma}

\begin{proof}
    The proof proceeds by induction on $\lambda$. The case $\lambda=0$ follows directly from the definition of $C_0$. Suppose that $S_{\lambda} \subseteq C_0(0)$ for some $\lambda \ge 0.$ For $\omega \in S_{\lambda}$, there exists some $W\in \mathcal{C}_0$ such that $\omega = W(0).$ From the relation $L\omega = [W,f_0](0)$, it follows that  $L\omega \in C_0(0)$. From the  $C_0(0)$-invariance of the distribution $C_0$ it follows that $[[f_0,\omega],W] \in \mathcal{C}_0$, thus $[[f_0,\omega],W](0) = 2\Phi(\omega)-DW(0)L\omega \in C_0(0)$. However, we note that $DW(0)L\omega \in C_0(0)$  since $[L\omega, W] \in \mathcal{C}_0$. Hence $\Phi(\omega) \in C_0(0)$. Therefore $\langle L\omega, \Phi(\omega) ~\rvert~ \omega \in S_{\lambda}\rangle \subseteq C_0(0)$. As a result,  $S_{\lambda+1} \subseteq C_0(0)$.
  \end{proof}

\begin{lemma}\label{C_0-Lemma2}
    Let $l \ge 0$. The iterated Lie bracket $B_l\coloneqq ad_{X_l}\cdots ad_{X_1}f_j$, where $X_i\in \{f_0,f_1,\cdots,f_{n-k}\}$ and $1\le j\le n-k$, is such that $B_l(0)\in S_{\lambda},~$ for some $\lambda \le l$.
\end{lemma}

\begin{proof}
    Here the proof proceeds by strong induction on $l.$ The case $l=0$ is vacuous since by convention $B_0=f_j \in S_0$ for any $1\le j\le n-k$. For $l=1$, we note that $B_1 = [X_1,f_j]$ for some $j$. Thus $B_1(0)=-Lf_j \in S_1 $ if $X_1=f_0$ and zero otherwise. 
    Now let $\kappa > 1$ and assume that $B_l(0) \in S_{\lambda}$ for some integer $\lambda \le l$, for all $0\le l <\kappa$.
    We have $B_{\kappa} = DB_{\kappa-1}X_{\kappa}-DX_{\kappa}B_{\kappa-1}$. Thus if $X_{\kappa}=f_0$, then $B_{\kappa}(0) = -LB_{\kappa-1}(0) \in S_{\lambda}$ for some $\lambda \le \kappa$ by the induction hypothesis. Otherwise, if $X_{\kappa}= f_j$ for $1\le j\le n-k$, then $B_{\kappa}(0) = DB_{\kappa-1}(0)f_{j} \in S_{\lambda}$ for some $\lambda \le \kappa-1$. To see why this last relation holds, we observe that $DB_{\kappa-1}f_{j} = D^2B_{\kappa-2}(X_{\kappa-1},f_{j})+DB_{\kappa-2}DX_{\kappa-1}f_{j}-DX_{\kappa-1}DB_{\kappa-2}f_{j}-D^2X_{\kappa-1}(B_{\kappa-2},f_{j})$. For $X_{\kappa-1}=f_0$, the  evaluation at the origin of the last term of this expression yields $ \Psi(B_{\kappa-2}(0),f_{j})$ which belongs to $S_{\lambda}$ for some $\lambda \le \kappa-1$ via the polarization identity in Proposition \ref{polarization-prop} since $B_{\kappa-2}(0) \in S_{\lambda}$ for some $\lambda \le \kappa-2$ by the induction hypothesis. Likewise, the first three terms involving the first and second derivative of the bracket $B_{\kappa-2}$ can be recursively reduced after $\kappa-3$ steps (for $\kappa >2$) to a linear combination of terms involving only the first derivative of $B_1$ since the vector fields $f_is$ are at most quadratic for $0\le i \le n-k$. Therefore $DB_{\kappa-1}(0)f_{j}$ indeed belongs to $S_{\lambda}$ for some $\lambda \le \kappa-1$. Hence the statement in the lemma follows by the principle of mathematical induction.
\end{proof}

Although the proof of the lemma above is relatively straightforward, the exact expression of the Lie bracket $B_l$, and that of its evaluation at the origin $B_l(0)$ becomes more involved the larger the order $l.$ We illustrate this by listing for instance in Table \ref{long-Brackets} and in the caption of Figure \ref{Bracket-Tree} the vectors $B_l(0)$ for the nonzero brackets $B_l, ~1\le l\le3$ shown in Figure \ref{Bracket-Tree} along with the $4$th order bracket $[f_l,[f_j,ad_{f_0}^2f_i]]$. These results can be obtained by direct computation from the definition of the Lie bracket or by repeated application of the Jacobi identity.

% \begin{table}[h!]
% \centering
%  \begin{tabular}{||c c c c c ||} 
%  \hline
%  $[f_0,f_i](0)$ & \textcolor{blue}{$[f_j,[f_0,f_i]]$} & $ad_{f_0}^2f_i(0)$ & $[f_0,[f_j,[f_0,f_i]]](0)$ & $ad_{f_0}^3f_i(0)$ \\ [0.5ex] 
%  \hline\hline
%  $-Lf_i$  & \textcolor{black}{$-\Psi(f_i,f_j)$}&  $L^2f_i$&$L\Psi(f_i,f_j)$ &$-L^3f_i$ \\ 
%  \hline
%  \end{tabular}
%  \caption{Value at the origin of the brackets in Figure \ref{short-brackets} except that of $[f_j,ad_{f_0}^2f_i]$ which is given in Table \ref{long-Brackets} since it's a bit more involved. Notice the general pattern of the $l$-th order ad-bracket $ad_{f_0}^lf_i(0)=(-1)^lL^lf_i$ which belongs to the subspace $S_l$.}
%  \label{short-brackets}
% \end{table}

\begin{remark}
    Figure \ref{Bracket-Tree} highlights the inherent graph structure of the accessibility algebra, More specifically, its generators make up a "forest" meaning a collection of rooted trees like the one shown in the figure. For instance, in the case of system $\Sigma_k$, the strong accessibility algebra  is made of $n-k$ such trees by Proposition \ref{Accessibility-Algebra}. The hope here in pointing this out is that it hints at how control scientists and engineers might leverage the powerful tools of graph learning in answering controllability questions.
    \end{remark}

 \begin{table}[h]
\begin{center}
\begin{tabular}{|c||c|}
\hline
$[f_j,ad_{f_0}^2f_i](0)$ & $L\Psi(f_i,f_j)+ \Psi(f_j,Lf_i)- \Psi(f_i,Lf_j)$ \\
\hline
$\mathbf{[f_l,[f_j,ad_{f_0}^2f_i]]}$ & $\Psi(\Psi(f_i,f_j),f_l)+ \Psi(\Psi(f_i,f_l),f_j)- \Psi(f_i,\Psi(f_j,f_l))$\\
\hline
\end{tabular}
\end{center}
\caption{Value at the origin of the  order $3$ bracket $[f_j,ad_{f_0}^2f_i]$ and  the order $4$  constant bracket $\mathbf{[f_l,[f_j,ad_{f_0}^2f_i]]}$. The corresponding vectors are explicitly  elements of $S_2$ which in turn is a subspace of $S_3$ and $S_4$ respectively by the polarization identity in Proposition \ref{polarization-prop}. Compared to the other brackets in Figure \ref{Bracket-Tree}, notice how the expressions here start to become more and more involved the higher the order of the bracket.}
\label{long-Brackets}
\end{table}

We define the set $\underset{\lambda \ge 0}{\sum}{S_{\lambda}}:= \big\{s_{\lambda_1}+\cdots+s_{\lambda_t} ~\rvert ~s_{\lambda_i}\in S_{\lambda_i},~ \lambda_i \ge 0,~t\ge1\big\}$ which is simply the collection of sums of elements of the $S_{\lambda}$, all of which are zero but finitely many. This set is endowed with a natural vector space structure induced by that of the $S_{\lambda}$. Furthermore, since the $S_{\lambda}$ are subspaces of $\mathbb{R}^n$, it follows that the increasing sequence $\big(S_{\lambda}\big)_{\lambda\ge 0}$ must  eventually be stationary. As a result, there must exist some integer $\kappa_0 \ge 0 $ such that $\underset{\lambda \ge 0}{\sum}{S_{\lambda}}=S_{\kappa_0}$. This motivates the following proposition.

\begin{proposition}\label{C_0-characterization}
    There exists some integer $0 \le \kappa_0 \le k$ such that $C_0(0)=S_{\kappa_0}.$
\end{proposition}

\begin{proof}
    The definition of $C_0$ and the combined  application of Proposition \ref{Accessibility-Algebra} followed by Lemma \ref{C_0-Lemma1} and  Lemma \ref{C_0-Lemma2} imply that $C_0(0)=\sum_{\lambda \ge 0}S_{\lambda}$. Furthermore, the sequence $S_{\lambda}$ is stationary from the first $\lambda$ such that $ \dim S_{\lambda} = \dim S_{\lambda+1}$. Hence, noting that $\dim S_{0} = n-k$, it follows that $\kappa_0 \le k$.
\end{proof}

We're now in the position of stating the main result of this section which is also the central result of this paper.
\begin{theorem}\label{Accessibitlity-Theorem}
    Let $1\le k \le n-1$. The control system $\Sigma_k$ is locally strongly accessible from the origin if and only if $S_k = \mathbb{R}^n.$
\end{theorem}

\begin{proof}
Fix $1\le k \le n-1$. We know from Lemma \ref{C_0-Lemma1} that $S_k \subseteq C_0(0)$. Thus $S_k = \mathbb{R}^n$ is sufficient for local strong accessibility of $\Sigma_k$ from the origin. Conversely, to show that this condition is also necessary, we'll first treat the case $k=1$ as baseline. So, we consider the change of coordinates $x = \delta_1 f_1+ \cdots +\delta_{n-1}f_{n-1} + \delta_nf_n$, where $f_n:=*\Big( \underset{j=1}{ \overset{n-1}{\wedge}}f_j\Big)$ is the Hodge-star of the $(n-1)$-vector obtained by wedging together the members of the family $\{f_1,\cdots,f_{n-1}\}$. We adopt the following notation :
           $ \Delta_a^{(1)} := a\odot\big(\delta_1P_2f_1+\cdots+\delta_{n-1}P_2f_{n-1}\big)\odot P_2f_n,~ 
            \Delta_a^{(2)} := a\odot P_2f_n \odot \big(\delta_1P_2f_1+\cdots+\delta_{n-1}P_2f_{n-1}\big),~
            \Delta_a^{(3)} := a\odot P_2f_n \odot P_2f_n, ~$ and $
            \Delta_a^{(4)} := a \odot \big(\delta_1P_2f_1+\cdots+\delta_{n-1}P_2f_{n-1}\big) \odot \big(\delta_1P_2f_1+\cdots+\delta_{n-1}P_2f_{n-1}\big)$.
            We also define $\Delta_b^{(j)}$ and $\Delta_c^{(j)}$ for $1\le j\le 4$ accordingly. With this notation and the change of coordinates above, the control system (\ref{class_S_Ctrl}) becomes
    \begin{align*}
        \begin{cases}
            \dot{\delta}_i = \frac{f_i^T}{||f_i||^2}\bigg(\underset{j=1}{\overset{n}{\sum}}\delta_j Lf_j +\underset{\alpha \in \{a,b,c\}}{\sum}\delta_n\big(\Delta_{\alpha}^{(1)}+\Delta_{\alpha}^{(2)}+\delta_n\Delta_{\alpha}^{(3)}\big)+\Delta_{\alpha}^{(4)} \bigg) + u_i~, \hspace{5pt} i=1,\cdots,n-1\\
            \vdots\\
            \dot{\delta}_n = \delta_n\frac{f_n^T}{||f_n||^2}\bigg(Lf_n+\underset{\alpha \in \{a,b,c\}}{\sum}\Delta_{\alpha}^{(1)}+\Delta_{\alpha}^{(2)}+\delta_n\Delta_{\alpha}^{(3)} \bigg)   + \frac{f_n^T}{||f_n||^2}\bigg(L\big( \underset{i=1}{\overset{n-1}{\sum}}\delta_if_i\big)+ \underset{\alpha \in \{a,b,c\}}{\sum} \Delta_{\alpha}^{(4)} \bigg) 
        \end{cases}
    \end{align*}
If the condition fails, then  the subspace $S_1 = Span\big\{f_1,\cdots,f_{n-1},L\omega +\Phi(\omega)~\big\rvert~ \omega\in S_0\big\}$  does not have maximum dimension. As a result, noting that $\Delta_a^{(4)}+\Delta_b^{(4)}+\Delta_c^{(4)} = \Phi(\sum_{j=1}^{n-1}\delta_jf_j)$ and using the mixed product property of the Hodge-star which yields a determinant in this case, it turns out that the second term in the $\dot{\delta}_n$ equation  vanishes identically. Therefore the transformed system cannot be accessible since  the hyperplane $\{\delta_n=0\}$ will be absorbing. It follows  from the invertibility of the change of coordinates that the control system $\Sigma_1$ is not accessible either. Now Let $1< k \le n-1$, and assume that $\Sigma_k$ is locally strongly accessible from the origin. There must exist some $w_1^0 \in S_1 \setminus S_0$; for if not, the sequence $(S_{\lambda})_{\lambda \ge 0}$ will be constant and equal to $S_0$ which in turn will imply that $C_0(0) = S_0$ by Proposition \ref{C_0-characterization}. This will contradict the accessibility assumption of $\Sigma_k$. So, the family $\{f_1,\cdots,f_{n-k},w_1^0 \}$ is linearly independent in $S_1$. If $S_1 = \mathbb{R}^n$ then we're done. Otherwise  $S_1 \neq\mathbb{R}^n$, and there must exist some $w_2^0 \in S_2 \setminus S_1$; for if not, the sequence $(S_{\lambda})_{\lambda \ge 1}$ will be constant and equal to $S_1$ implying that $C_0(0) = S_1$ which again will contradict the assumption of accessibility of $\Sigma_k$. So, we have the linearly independent family $\{f_1,\cdots,f_{n-k},w_1^0,w_2^0 \}$ in $S_2.$ This procedure can be continued to exhibit vectors $w_j^0 \in S_j \setminus S_{j-1}$ for $3 \le j \le k-1$ such that the family $\{f_1,\cdots,f_{n-k},w_1^0,\cdots w_{k-1}^0 \}$ is linearly independent in $S_{k-1}.$ From here, we  proceed as in the base case $k=1$.
\end{proof}

The proof above shows that the  condition in the theorem is  necessary for accessibility in the wider sense namely (AP), thus we have the following proposition.
\begin{proposition}\label{acc-Proposition-NC}
    The control system $\Sigma_k$ has the accessibility property from the origin only if $S_k = \mathbb{R}^n~$, $1\le k \le n-1$. 
\end{proposition}

\begin{example}
   Consider the affine quadratic control system in $\mathbb{R}^5$ afforded by $L = diag(0,1,-1,0,0)+\mathbf{E}_{45},~$the vectors $a=\mathbf{0},~b=(0,0,0,-1,-1)^T,~ c=(0,0,1,0,1)^T,~$ and the control vector $f_1=\frac{\partial}{\partial x_1}$. This system can also be written as:
   $$\Sigma_4:\quad\begin{cases}
        \dot{x}_1=u_1\\\dot{x}_2=x_2\\\dot{x}_3=x_4x_5-x_3\\
        \dot{x}_4=x_5-x_1^2\\\dot{x}_5=x_1x_2-x_2^2
    \end{cases}$$
    A straightforward computation shows that $S_4 = \langle\mathbf{e}_1,\mathbf{e}_4\rangle$, so this system is not accessible from the origin. 
\end{example}

\begin{remark}
    One of the main implications of the result in Theorem \ref{Accessibitlity-Theorem} is that it quantifies  the \emph{"degree of reachability"} of the system $\Sigma_k$ via the dimension of the subspace  $S_k$. That is to say $dim(S_k)$ measures the dimension of the largest neighborhood the system can explore. For instance, in the example above the system can at most explore a $2$-dimensional neighborhood, therefore more actuators are needed to explore additional dimensions of the state space.
\end{remark}

\subsection{Application: The Sprott System}

Here we highlight a corollary of the above theorem to the subclass of systems where $a=\mathbf{1}$ and $b=\mathbf{0}=c$ to which the Sprott system (\ref{eq:SprottSystem-Delta-dynamics}) belongs. We begin with a couple of useful lemmas.

\begin{lemma}\label{lemm:lemmaUseful_1}
    Let $u,v$ be non-colinear vectors in $\mathbb{R}^3$. There exists  $\omega \in Span\{ u,v\}$ such that the set of vectors $\mathcal{S} = \{u,v,\Delta(P_2\omega)P_2\omega\}$ is linearly independent.
\end{lemma}
\begin{proof}
     Let $\omega = \alpha u +\beta v$ for some $\alpha,\beta \in \mathbb{R}$, and let $u\times v = (d_1,d_2,d_3)^T$. The determinants $d_i$ are not all vanishing. So, for simplicity we'll assume that $d_i \neq 0$, for  $1\le i \le 3.$ We note that $det(\mathcal{S}) = d_1(\alpha u_2+\beta v_2)^2+d_2(\alpha u_3+\beta v_3)^2 + d_3(\alpha u_1+\beta v_1)^2$. Suppose for the sake of contradiction that the assertion is false, namely $det(\mathcal{S}) =0$ for all $\alpha,\beta \in \mathbb{R}$. The proof proceeds by distinguishing the three main cases afforded by the signs of the $d_i$. First, the $d_i$ are all of the same sign. We can assume for instance that $d_i >0$ for $1\le i \le 3$. $det(\mathcal{S}) = 0$ implies that $d_1(\alpha u_2+\beta v_2)^2 = 0$. In particular, for $\alpha = v_3$ and $\beta = -u_3$, we arrive at $d_1=0$ which is a contradiction. Next, exactly one of the $d_i$ is negative. We can assume for instance that $d_1<0$ and $d_2,d_3 >0$. Setting $\alpha = v_2$ and $\beta = -u_2$, we get $d_2d_1^2+d_3^3=0$ which yields again a contradiction. Finally, exactly two of the $d_i$ are negative. We can assume for instance that $d_1,d_2<0$ and $d_3 >0$. Here also, setting $\alpha = v_1$ and $\beta = -u_1$ we arrive at $d_1d_3^2+d_2^3 = 0$ leading to yet again a contradiction. In all these exhaustive cases we reach a contradiction. Hence it must be the case that the assertion in the Lemma holds. 
\end{proof}

\begin{lemma}\label{lemm:lemmaUseful_2}
    The vectors $w$ and $\Delta(P_2w)P_2(w)$ in $\mathbb{R}^n$ are linearly dependent  if and only if $w \in Span\{\mathbf{1}\}$. 
\end{lemma}
\begin{proof}
 The subspace $Span\{\mathbf{1}\}$ is invariant under the map $w \mapsto \Delta(P_2w)P_2w$. Conversely, suppose $\exists \lambda \in \mathbb{R}$ such that $w = \lambda \Delta(P_2w)P_2w$. If $\lambda = 0$ we're done. Otherwise, $\lambda \ne 0$, and the following holds: $w_i = \lambda^{2^n-1}w_i^{2^n}$ for $1\le i \le n$. This implies that $w_i=0$ or $w_i = \frac{sign(\lambda)}{\lambda}$ for all $i$. Thus $ w \in Span\{\mathbf{1}\}$ in either case. 
\end{proof}

 \begin{proposition}\label{prop:propositionSubClass_Sprott}
     Let $a = \mathbf{1}$,~ $b= \mathbf{0}= c$ in $\mathbb{R}^3$ and assume that the diagonal $Span\{\mathbf{1}\}$ is $L$-invariant. The control system $\Sigma_k$ is locally strongly accessible from the origin if and only if
    \begin{enumerate}
        \item $f_1$ and $f_2$ are linearly independent for the case $k=1$.
        \item $f_1 \not\in Span\{\mathbf{1}\}$ for the case $k=2$.
    \end{enumerate}
\end{proposition}
\begin{proof}
    $Span\{ u,v,L\omega + \Delta(P_2\omega)P_2\omega ~\rvert~ \omega \in Span\{ u,v \} \} = \mathbb{R}^3$  if and only if $u$ and $v$ are linearly independent by Lemma \ref{lemm:lemmaUseful_1}, and the  first assertion follows from Theorem \ref{Accessibitlity-Theorem}. To prove the second assertion, we note that if the diagonal $Span\{\mathbf{1}\}$ is $L$-invariant, then Lemmas \ref{lemm:lemmaUseful_1} and \ref{lemm:lemmaUseful_2} imply in particular that $f_1 \notin Span\{1\}$ if and only if  $Span\{ f_1,\Delta(P_2f_1)P_2f_1,L\omega + \Delta(P_2\omega)P_2\omega ~\rvert~ \omega \in \langle f_1,\Delta(P_2f_1)P_2f_1 \rangle \} = \mathbb{R}^3$. Conclude with Theorem \ref{Accessibitlity-Theorem}.
\end{proof}

\begin{proposition}
    The Sprott system (\ref{eq:SprottSystem-Delta-dynamics}) is locally strongly accessible from the origin with two control vectors (resp: one control vector) if and only if they are linearly independent (resp: it is not along the diagonal).
\end{proposition}

\subsection{The Rigid Body and Crouch's Results Consistency}\label{subsec:Crouch-consistency}

We've shown in section \ref{subsection:some_Examples} that class $\mathcal{C}$ in (\ref{class_S}) contains the rigid body angular velocity dynamic, and we gave in Theorem \ref{Accessibitlity-Theorem} necessary and sufficient conditions for strong accessibility for systems in this class. However Crouch \cite{crouch1984spacecraft}  gave necessary and sufficient conditions for  accessibility of the rigid body dynamics which he coupled to the Poisson stability of the drift vector field to conclude controllability. So, a natural question is that of whether or not our condition recovers Crouch's in the case of the rigid body. The answer is yes, and we present the details below by carefully relating the rigid body model \cite{marsden2013introduction}\cite{schaub2003analytical} to our formalism. 
The rigid body  angular velocity dynamics is given by the equation 
$I\dot{\omega}
        = S(\omega)I\omega +\tau,$ 
        where $S(\omega)=\hat{\omega}$ is the skew symmetric matrix given in section \ref{par:RB} above. Since $I$ is invertible, in fact it's symmetric positive definite, we can equivalently write this dynamic as:  $\dot{\omega} = I^{-1}S(\omega)I\omega + I^{-1}\tau$. The main implication of this equation  is that  whenever we consider the control affine form of the rigid body model, the control vector terms  represents the \emph{scaled} total torques applied to the body. We know from Crouch \cite{crouch1984spacecraft} that $ Span \big\{ b_1,b_2,S(\omega)J^{-1}\omega,~ \omega \in span\{b_1,b_2\} \big\} = \mathbb{R}^3$ is necessary and sufficient for local strong accessibility of the rigid body dynamic. 
        
        We now show that this condition is equivalent to the one given in Theorem \ref{Accessibitlity-Theorem} for $k=1$. We recall from (\ref{par:RB})  that for the rigid body, $\Phi(x) = \Delta_{\frac{1}{\xi}}S(x)\Delta_{\xi}x$ for some $\xi = (\xi_1,\xi_2,\xi_3)^T \in \mathbb{R}^3$ with $\xi_i > 0$, where  $\Delta_{\xi}$ is the principal moment of inertial matrix of the rigid body with $\Delta_{\frac{1}{\xi}} = \Delta_{\xi}^{-1}$. From the above discussion, it's now clear that the control term $u_if_i$ has the physical meaning of scaled torque, namely $f_i = \Delta_{\frac{1}{\xi}}b_i$ for some $b_i$, $1\le i\le2$. Therefore it follows that for any $\alpha, \beta \in \mathbb{R}$: $\Delta_{\xi}\Phi(\alpha f_1+\beta f_2) = \Delta_{\xi}\Delta_{\frac{1}{\xi}}S(\alpha f_1+\beta f_2)\Delta_{\xi}(\alpha f_1+\beta f_2)
        = -\Delta_{\frac{1}{\xi}}(\alpha b_1+\beta b_2)\times (\alpha b_1+\beta b_2)
        = -S(\alpha b_1+\beta b_2)\Delta_{\frac{1}{\xi}}(\alpha b_1+\beta b_2)$. It results from the invertibility of $\Delta_{\xi}$ that $Span \big\{ f_1,f_2,\Phi(\alpha f_1+\beta f_2)~\rvert~\alpha,\beta \in \mathbb{R} \big\} = \mathbb{R}^3$ if and only if $Span \big\{ b_1,b_2,S(\alpha b_1+\beta b_2)\Delta_{\frac{1}{\xi}}(\alpha b_1+\beta b_2)~\rvert~\alpha,\beta \in \mathbb{R} \big\} = \mathbb{R}^3$ which is precisely the Crouch condition in the case of two pairs of gas jets externally applied to the body. The argument is essentially the same in the case $k=2$ (i.e a single pair of gas jets) by making the substitution $f_2 = \Phi(f_1)$ accordingly.

\subsection{Reduction to Linear Systems}\label{subsec:Kalman-consistency}

If we turn off the quadratic component $\Phi(x)$ of the drift in $\mathcal{C}$ (\ref{class_S}), meaning if $\mathbf{a}=\mathbf{0}=\mathbf{b}=\mathbf{c}$, then class $\mathcal{C}$ reduces to that of linear systems, and system (\ref{class_S_Ctrl}) becomes  $\Sigma_k:  \dot{x} = Lx+B_ku$, where $B_k$ is the matrix whose columns are the linearly independent constant  control vectors $\{f_{i}\}_{1\le i\le n-k}$, for any $1\le k \le n-1$.

On the one hand, the Kalman rank condition gives a necessary and sufficient condition for the controllability of $\Sigma_k$. On the other hand, in light of Proposition \ref{CtrbEquvLocStrAccForLinSys} which establishes the equivalence between strong accessibility and controllability for linear systems, Theorem \ref{Accessibitlity-Theorem} also gives a necessary and sufficient condition for controllability of $\Sigma_k$. The natural question that arises is that of whether the condition $S_k=\mathbb{R}^n$ given in Theorem \ref{Accessibitlity-Theorem} recovers the Kalman rank condition for linear systems. The answer is yes, and we give the details below.

We recall the following propositions. The first one is a fairly classic result that builds on the work of Sussmann and Jurdjevic  \cite{sussmann1972controllability}, and Krener \cite{krener1974generalization}. Meanwhile, the second one establishes the equivalence between the notions of controllability and strong accessibility for linear systems.

\begin{proposition}[Proposition 2.1 in \cite{sussmann1987general}]\label{SussmannLARCnAP}
    Let $\mathcal{F}$ be a family of $C^{\infty}$ vector fields on a $C^{\infty}$ manifold M. Then  the \emph{LARC} at p implies the \emph{AP} from p. Conversely, the \emph{AP} from p implies the \emph{LARC} at p if M is a real-analytic manifold and the members of $\mathcal{F}$ are real-analytic
\end{proposition}

\begin{proposition}\label{CtrbEquvLocStrAccForLinSys}
    For linear systems $\dot{x}=Lx+Bu$ in $\mathbb{R}^n$, the notions of controllability and local strong accessibility from the origin are equivalent.
\end{proposition}

\begin{proof}
   First, note that the strong accessibility distribution at the origin is given by the reachable subspace $Span\{B,LB,\cdots,L^{n-1}B\}$ via Cayley-Hamilton. Hence the controllability is equivalent to LARC at the origin. The same argument combined with Proposition \ref{LocStrgAccSuffCond} also shows that  controllability implies local strong  accessibility from the origin. Next, it follows from Proposition \ref{SussmannLARCnAP} that LARC and AP are equivalent at the origin since the manifold $\mathbb{R}^n$ and the vector fields involved are real-analytic. Finally, local strong accessibility from the origin implies AP from the origin by definition.
   \end{proof}
   
\begin{figure}[H]
    \centering
    \includegraphics[width=.65\textwidth]{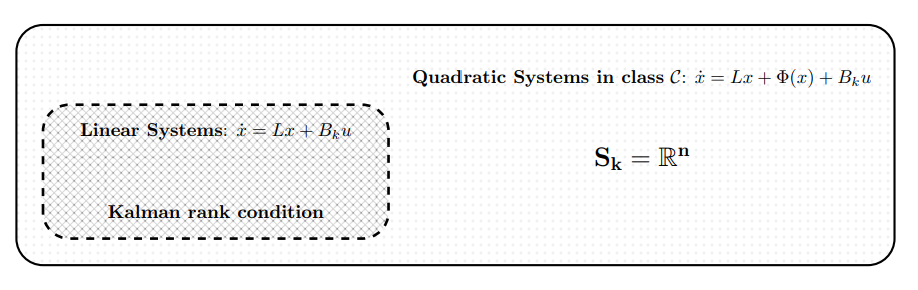} % Image using textwidth to define size
    \caption{ The Kalman rank condition covers the controllability of the class of linear systems which is a subclass of the class of quadratic systems $\mathcal{C}$. $S_k=\mathbb{R}^n$ covers the strong accessibility of the class $\mathcal{C}$, and it recovers the Kalman rank condition in the sense that when crossing the dashed boundary into linear systems, this condition  is equivalent to the Kalman rank condition. Therefore, the condition  $S_k=\mathbb{R}^n$ extends the Kalman rank condition to $\mathcal{C}$ for $1\le k \le n-1$.} 
    \label{KCC-generalized}
 \end{figure}

 We refer the interested  reader to classic reference by Kawsky \cite{kawski2017high} for a comprehensive survey on notions of controllability, and \cite{boscain2023local}  for a more recent work on the same topic. The following theorem connects the classic Kalman controllability result to the main result of this paper, namely Theorem \ref{Accessibitlity-Theorem}. See Figure \ref{KCC-generalized} for a pictorial representation.

\begin{theorem}\label{KCCequivWkeqRn}
    For $\Phi(x)=0$, consider the corresponding linear control system  $\Sigma_k: \dot{x}=Lx+B_ku$ defined in  $\mathbb{R}^n$. Then for $1\le k\le n-1$,  the Kalman rank condition  $Span\{B_k,LB_k,\cdots,L^{n-1}B_k\} = \mathbb{R}^n$ is equivalent to the condition  $S_k=\mathbb{R}^n$. 
\end{theorem}

\begin{proof}
    Let $1\le k\le n-1$, and let's denote by (K) the Kalman rank condition. Note that for linear systems (i.e $\Phi(x)=0$), $S_k=Span\{B_k,LB_k,\cdots,L^kB_k\}$ where the subspace $S_k$ was defined in (\ref{FundamentalSubspaces}). Hence it's straightforward to see that $S_k=\mathbb{R}^n$ is sufficient for (K). Conversely, to show that $S_k=\mathbb{R}^n$ is necessary for (K), we proceed inductively on $k$. The case $k=n-1$ is clearly vacuous. For the case $k=1$, given that the dimension of the subspace $Span\{B_1\}$ is not maximum, (K) implies that $Span\{LB_1\} \not\subset Span\{B_1\}$ which in turn implies that $S_1=\mathbb{R}^n$ since $Span\{B_1\}$ is of codimension $1$. For the case $k \ge 2$, let $1 \le j< k$ and proceed inductively on $j$ as follows: Define the subspace $V_{j,k} = Span\{B_k,LB_k,\cdots,L^jB_k\}$. At step 1, given that $rank(B_k)=n-k < n$, (K) implies that $\dim V_{1,k}\ge (n-k)+1$. At step 2, if $Span\{L^2B_k\} \subseteq V_{1,k}$, then $Span\{L^iB_k\} \subseteq V_{1,k}$ for all $2\le i\le n-1$. Thus $Span\{B_k,LB_k,\cdots,L^{n-1}B_k\} \subseteq V_{1,k} \subset S_k$. Hence (K) implies $S_k=\mathbb{R}^n$. Otherwise, $Span\{L^2B_k\} \not\subset V_{1,k}$, and we must have $\dim V_{2,k} \ge (n-k)+2$. The argument in step 2 can be eventually repeated accordingly until the subspace $V_{k-1,k}$ satisfying $\dim V_{k-1,k} \ge (n-k)+(k-1)$ is constructed at step k-1. Finally, if $Span\{L^kB_k\} \subseteq V_{k-1,k}$ then the result follows by the same argument as in step 2. Otherwise, we have $Span\{L^kB_k\} \not\subset V_{k-1,k}$ which implies that $\dim S_k \ge (n-k)+(k-1)+1 = n$. Hence $S_k=\mathbb{R}^n$.
\end{proof}

\begin{remark}
    The equivalence in Theorem \ref{KCCequivWkeqRn} can alternatively be ascertained through a sequence of pre-established results as follows: First, the Kalman rank condition is equivalent to LARC at the origin. The core argument of this assertion was already given in the proof of Proposition \ref{CtrbEquvLocStrAccForLinSys}. See also \cite{bloch2015nonholonomic} and the discussion in chapter 4 for instance for further details. Second, LARC at the origin is equivalent to AP from the origin for $\Sigma_k$, $1\le k\le n-1$. This is essentially Proposition \ref{SussmannLARCnAP}. Third, AP from the origin implies $S_k = \mathbb{R}^n$. This is the content of Proposition \ref{acc-Proposition-NC} which can be seen as a  corollary of the central theorem of this paper, namely Theorem \ref{Accessibitlity-Theorem}. Finally, the loop closes itself trivially by noticing that for linear systems $\Sigma_k$, we have $S_k=Span\{B_k,LB_k,\cdots,L^kB_k\}$. Hence we have indeed $S_k=\mathbb{R}^n$ implies the Kalman rank condition.
\end{remark}

\section{Small Time Local Controllability}

In this section, we analyze the Small-Time Local Controllability (STLC) of the control affine system $\Sigma_k$. More specifically, we're interested in what can be said about its controllability under the condition $S_k = \mathbb{R}^n$ which as we know from Theorem \ref{Accessibitlity-Theorem} is necessary and sufficient for strong accessibility. First, we show through  a simple counterexample that this condition is not in general sufficient for local controllability. However as a  second step, we show that under some relatively mild assumptions on either the  linear part $L$ or the homogeneous quadratic part $\Phi$ of the drift, this condition becomes necessary and sufficient for small time local controllability at least in the case $k=1$ which we refer to as the \emph{rank $1$ underactuation} case. Finally, we show that in the case $k=n-1$, namely for single input systems, if the linear part $L$ of the drift vanishes identically then the system cannot be STLC. 

Let's start by recalling some useful background results about local controllability of control systems. More detailed  information can be found in \cite{bloch2015nonholonomic} \cite{nijmeijer1990nonlinear} \cite{coron2007control} \cite{bullo2019geometric} and the references therein. 
Currently we're not aware of any necessary and sufficient condition to check for small-time local controllability for general control systems. However,  there are some useful conditions that are either necessary or sufficient, and we present some of them below. We start with sufficient condition theorems. The first one guarantees STLC of a nonlinear control system from the controllability of its linearization.

\begin{theorem}[Theorem 3 in  \cite{markus1965controllability}]\label{theorem:STLC_from_linearization}
Let $(x_e,u_e)$ be an equilibrium point of a nonlinear control system $\dot{x}=f(x,u)$. Let's assume that its linearization at $(x_e,u_e)$ is controllable. Then the nonlinear control system is small-time locally controllable from $(x_e,u_e)$. 
\end{theorem}

The next theorem (\ref{theorem:Sussmann_sufficiency}) is a powerful general sufficient result due to Sussmann \cite{sussmann1987general} which extends some previous work on the subject matter among which that of Hermes \cite{hermes1976local}\cite{hermes1978lie} \cite{hermes1982local}\cite{hermes1982control}, Crouch and Byrnes \cite{crouch1986symmetries},  Brunovsky \cite{brunovsky1974local}, Jurdjevic \cite{jurdjevic1983polynomial},  and many others. See for instance Bloch et al \cite{bloch1992Reyanuglu} for a nice application of this result to nonholonomic dynamic systems.

Let $\mathbf{X}= (X_0,\cdots,X_m )$ be a finite sequence of indeterminate, and denote by $A(\mathbf{X})$ the free real associative  algebra generated by the $X_i$s. Let $L(\mathbf{X})$ denote the Lie subalgebra of $A(\mathbf{X})$ generated by $X_0,\cdots, X_m$. We denote by $Br(\mathbf{X})$  the smallest subset of $L(\mathbf{X})$ that contains $X_0,\cdots, X_m$ and is closed under bracketing. We refer as \emph{bracket} an element of $Br(\mathbf{X})$, namely an element of $L(\mathbf{X})$  that cannot be written as a linear combination of other elements. Given a bracket $B$, we define the \emph{degree} of $B$ to be $\delta(B) = \sum_{i=0}^{m}\delta^i(B)$, where $\delta^i(B)$ denotes the number of times the indeterminate $X_i$ occurs in the bracket $B$. The bracket $B$ is said to be \emph{bad} if $\delta^0(B)$ is odd meanwhile $\delta^i(B)$ is even for all $i=1,\cdots, m$. The bracket is said to be \emph{good} if it's not bad. Let $\theta \in [1,\infty]$, we define the $\theta$-degree of the bracket $B$ as follows:
\begin{align}
\delta_{\theta}(B)=\begin{cases}
\frac{1}{\theta} \delta^0(B)+\sum_{i=1}^m \delta^i(B), & \theta \in [1,\infty),\\
\sum_{i=1}^m \delta^i(B), & \theta = \infty.
\end{cases}
\end{align}
Consider the $(m+1)$-tuple $\mathbf{f} = \{f_0,\cdots,f_m\}$ of $C^{\infty}$ vector fields on a $C^{\infty}$ manifold $M$. Each element $f_i$ is a member of $D(M)$ the algebra of all partial differential operators on the space $C^{\infty}(M)$ of real-valued functions on $M$. We define the evaluation map 
$Ev(\mathbf{f})\colon A(\mathbf{X}) \longrightarrow D(M)$ obtained by substituting the $f_i$ for the $X_i$, namely
   $Ev(\mathbf{f}) \big(\sum_I{a_IX_I} \big) = \sum_I{a_If_I}$, 
where $f_I = f_{i_1}f_{i_2}\cdots f_{i_s}$ for $I = (i_1,\cdots,i_s)$. Finally, for $x_0 \in M$, we denote by 
    $Ev(\mathbf{f})_{x_0}\colon L(\mathbf{X}) \longrightarrow T_{x_0}M$
the map that evaluates at $x_0$ the vector field $Ev(\mathbf{f})(B)$. We have the following definition

\begin{definition}\label{def:Suss_Neutralization}
Let $\theta \in [1,\infty]$. The control-affine system $\dot{x} = f_0(x)+\sum_{i=1}^m u_if_i(x)$ satisfies the Sussmann condition $S(\theta)$ if for every bad bracket $B$ there exist brackets $C_1,\cdots,C_k$ such that 

$i)~ Ev(\mathbf{f})_{x_0}(\beta(B)) \in Span\big\{Ev(\mathbf{f})_{x_0}(C_j),~ j=1,\cdots,l \big\}$

$ii)~ \delta_{\theta}(C_j) < \delta_{\theta}(B)$ for all $j=1,\cdots, l$.
\end{definition}
where  $\beta(B) = \underset{\sigma \in S_m}{\sum}\sigma(B)$ is the symmetrization  of the bracket $B$.

\begin{theorem}[Theorem 7.3 in \cite{sussmann1987general}]\label{theorem:Sussmann_sufficiency}
Consider the control-affine system $\dot{x} = f_0(x)+\sum_{i=1}^m u_if_i(x)$, and let $x_0$ be an equilibrium for the drift term $f_0$. Assume that $\mathbf{f} = \{f_0,\cdots,f_m\}$ satisfies the Lie Algebra Rank Condition at $x_0$, and further assume that there exists $\theta \in [1,\infty]$ such that the system satisfies the Sussmann condition $S(\theta)$. Then the system is small-time locally controllable from $x_0$.
\end{theorem}

We close this brief background survey by recalling a couple of very useful necessary condition results for STLC that will be needed subsequently.

\begin{theorem}[\cite{hermann1963accessibility}\cite{nagano1966linear}\cite{sussmann1973orbits}]\label{theorem:NC_STLC_analytic}
Assume that the control system $\dot{x}=f(x,u)$ is small-time locally controllable at the equilibrium $(x_e,u_e)$ and further assume that $f$ is analytic. Then the control system  satisfies the Lie Algebra Rank Condition \emph{(LARC)} at $(x_e,u_e)$.
\end{theorem}

This necessary condition is also sufficient for  linear control systems \cite{kalman155070controllability} and  drifless affine nonlinear control  systems \cite{rashevsky1938any}\cite{chow1939uber}.

The following result due to Sussmann \cite{sussmann1983lie}  which builds on the work of Hermes \cite{hermes1977controlled} will be referred  here as the \emph{Hermes-Sussmann necessary condition}. It's also known as the Clebsh-Legendre condition in optimal control \cite{kawski2017high}. This result is particularly useful in dealing with single-input affine control systems.
Let Lie$(f_0,f_1)$ be the Lie algebra generated by the smooth vector fields $f_0$ and $f_1$. For every integer $j\ge 0$, we define the linear subspace $\mathcal{S}^j(f_0,f_1)$ of Lie$(f_0,f_1)$ as the span of all Lie monomials in $f_0$ and $f_1$ which involve $f_1$ at most $j$ times. For instance, we have $\mathcal{S}^1(f_0,f_1) = Span\{f_0, ad_{f_0}^kf_1,~ k\ge 0 \}$.

\begin{theorem}[Theorem 6.3 in \cite{sussmann1983lie}]\label{theorem:NC_single-input-Sussmann}
    Consider an analytic single-input affine control system $\dot{x} = f_0(x)+ uf_1(x)$, and consider a point $x_0$ on the state-space for which  $[f_1,[f_0,f_1]](x_0) \notin \mathcal{S}^1(f_0+\bar{u}f_1,f_1)(x_0) $. Then this control system is not small-time locally controllable from $x_0$.
\end{theorem}

% We now illustrate through a counterexample that the condition $W_k=\mathbb{R}^n$ is not in general sufficient for STLC. So, 
\begin{counterexample}
    Consider the control system $\Sigma_1$ in $\mathbb{R}^3$ afforded by $L=0$, $b=(0,0,1)^T$ and $a=\mathbf{0}=c$ with control vectors fields $f_1=\frac{\partial}{\partial x}$ and $f_2=\frac{\partial}{\partial y}$. Equivalently this system reads: $$\Sigma_1:\quad\begin{cases}
        \dot{x}=u_1\\\dot{y}=u_2\\\dot{z}=y^2\\
    \end{cases}$$
We have $\langle f_1,f_2,\Phi(f_2)\rangle = \mathbb{R}^3$ since $\Phi(f_2) = \frac{\partial}{\partial z}$, hence $S_1=\mathbb{R}^3$ which by Theorem \ref{Accessibitlity-Theorem} ascertains the accessibility of $\Sigma_1$ from the origin. However, it's also straightforward to see that $\Sigma_1$ is not controllable since the $z$-component of the state variable is non-decreasing. Hence the condition $S_k=\mathbb{R}^n$ is not in general sufficient for STLC.

\end{counterexample}

The following theorem which is the main result of this section shows that under some additional assumptions on the linear or the homogeneous quadratic component of the drift, the condition in Theorem \ref{Accessibitlity-Theorem} also guarantees small-time local controllability.

\begin{theorem}\label{STLC-Theorem}
    Consider the control system $\Sigma_1$, and assume either that the subspace $S_0$ is not $L$-invariant  or $\Phi(f_i) \in S_0$ for all $1\le i \le n-1$, then $\Sigma_1$ is small-time locally controllable  from the origin if and only if $S_1=\mathbb{R}^n$.
\end{theorem}

\begin{proof}
    It follows from Proposition \ref{acc-Proposition-NC} that the condition $S_1=\mathbb{R}^n$ is  necessary for accessibility at the origin in the wider sense. Therefore it's also necessary for STLC from the origin by Theorem \ref{theorem:NC_STLC_analytic} and 
 Proposition \ref{SussmannLARCnAP}. Conversely, we  recall that the condition $S_1=\mathbb{R}^n$ is sufficient for LARC by Theorem \ref{Accessibitlity-Theorem}. Now, consider the same change of coordinates introduced in the proof of Theorem \ref{Accessibitlity-Theorem}. Denoting by $F_{\kappa}$, $0\le \kappa \le n-1$ the resulting vector fields, a direct computation yields the following brackets: $[F_j,F_0](0) = \sum_{i=1}^n{\frac{f_i^TLf_j}{||f_i||^2}}\mathbf{e}_i$,  $~[[F_0,F_j],F_l] = \sum_{i=1}^n{\frac{f_i^T\Psi(f_j,f_l)}{||f_i||^2}}\frac{\partial}{\partial_i}$ for $1\le j,l\le n-1$, and in  particular for $j=l$ we get $[[F_0,F_j],F_j]= \sum_{i=1}^n{\frac{2f_i^T\Phi(f_j)}{||f_i||^2}\frac{\partial}{\partial_i}}$. If $S_0$ is not $L$-invariant, then  $S_0 \cup \{[f_{j_0},f_0](0)\}$ neutralizes any potential bad bracket at the origin for some $j_0 \in \{1,\cdots,n-1\}$ resulting in the STLC of $\Sigma_1$ by Theorem \ref{theorem:Sussmann_sufficiency} with $\theta \in [1,\infty]$. However if $S_0$ is $L$-invariant, then the hypothesis imposes that $f_n^T\Phi(f_j)=0$ for all $1 \le j \le n-1$, implying that all type (1,2) bad brackets, namely brackets with one copy of the drift and two copies of the same control vector are neutralized by $S_0$. Furthermore, there must exist $j_0 \neq l_0$ such that $f_n^T\Psi(f_{j_0},f_{l_0}) \neq 0$ otherwise $S_0$ will become $\Phi$-invariant. To see this, it suffices to expand $\Phi(\omega) = \frac{1}{2}\Psi(\omega,\omega)$ for any $\omega \in S_0$ leveraging the symmetric bilinearity of the map $\Psi$. Nevertheless, the combined $(\Phi$ and $L)$-invariance of $S_0$ will contradict the premise $S_1=\mathbb{R}^n$. Therefore, any potential higher order bad bracket is neutralized by  $S_0\cup \{[f_{l_0},[f_{j_0},f_0]]\}$. Hence, we conclude by Theorem \ref{theorem:Sussmann_sufficiency} with $\theta \in [1,\infty)$ that  $\Sigma_1$  is STLC from the origin. 
\end{proof}

\begin{example}
    Consider the system $\Sigma_1$ in $\mathbb{R}^3$ given by $L = diag(2,0,-1),~$$\mathbf{a}=(0,0,1)^T,~\mathbf{b}=(0,1,1)^T,~\mathbf{c}=(-1,-2,0)^T,~ f_1=\frac{\partial}{\partial y}~$ and  $f_2=\frac{\partial}{\partial z}$. Equivalently, this system takes the form
    $$\Sigma_1:\quad\begin{cases}
        \dot{x}=2x-yz\\\dot{y}=x^2-2xz+u_1\\\dot{z}=y^2+x^2-z+u_2\\
    \end{cases}$$
    We note that $LS_0=\langle \mathbf{e}_3\rangle$, thus the first assumption in Theorem \ref{STLC-Theorem} is not met and the system is not linearly controllable. However $\Phi(f_1)=\mathbf{e}_3$ and $\Phi(f_2)=\mathbf{0}$, thus the second assumption in the theorem is met instead. Furthermore, we have $\langle f_1,f_2,\Phi(f_1+f_2)=\mathbb{R}^3$, thus $S_1=\mathbb{R}^3$ and the system is STLC from the origin by Theorem \ref{STLC-Theorem}.
\end{example}

\begin{proposition}\label{prop:Not-STLC-single-input-L=0}
    The affine control system $\Sigma_{n-1}$ with $L= 0$ in the drift  is not small-time locally controllable from the origin.
\end{proposition}
\begin{proof}
      The subspace $\mathcal{S}^1(f_0,f_1)(0) = Span\big\{f_0(0), ad_{f_0}^kf_1(0),~ k\ge 0 \big\}$ reduces to the line $S_0=Span\big\{ f_1\big\}$ under the assumption $L=0$.  If $ad_{f_1}^2f_0 (0) \notin \mathcal{S}^1(f_0,f_1)(0)$, then the result follows from Theorem \ref{theorem:NC_single-input-Sussmann} since the vector fields  $f_0$ and $f_1$ are  analytic. Otherwise, the subspace $S_{n-1}$ does not have maximum dimension since $\Phi(f_1)=\frac{1}{2}ad_{f_1}^2f_0(0)\in S_0$. Therefore the system is not even accessible  by Theorem \ref{Accessibitlity-Theorem}.
\end{proof}

The two cases considered in the proof of proposition above are not superfluous. In fact, they're complementary as the following two examples illustrate.
\begin{example}
    Consider the system $\dot{x} = \Delta_c(P_2x)P_3x+ uf_1$, where $c = (c_1,c_2,c_3)^T $ with $c_i > 0$ for $1\le i\le 3$. Define $\alpha = \sqrt{c_2c_3},~\beta = \sqrt{c_3c_1}$, $\gamma = \sqrt{c_1c_2}$ and the family $\mathcal{F}= \Big\{\big(\alpha^{-1},\beta^{-1},\gamma^{-1}\big),\big(-\alpha^{-1},-\beta^{-1},\gamma^{-1}\big),\big(-\alpha^{-1},\beta^{-1},-\gamma^{-1}\big),\big(\alpha^{-1},-\beta^{-1},-\gamma^{-1}\big) \Big\}$. For any control vector $f_1 \in \mathcal{F}$,
    we have  $ad_{f_1}^2f_0 \in \mathcal{S}^1(f_0,f_1)(0)$, namely the Hermes-Sussmann necessary condition is satisfied. An example of such system is the one afforded by $c = (1,1,1)^T$, namely 
    $\{ \dot{x}=yz,\dot{y}=xz,\dot{z}=xy\}$.
   This system was studied in \textnormal{\cite{chen2021controllability}} in the context of the controllability of hypergraphs. Proposition \ref{prop:Not-STLC-single-input-L=0} tells us that this system is not single-input STLC from the origin. How about accessibility? In fact, it's not too hard to see from Theorem \ref{Accessibitlity-Theorem} that this system is single-input locally strongly accessible from the origin if and only if the actuator $f_1=(f_{1_1},f_{1_2},f_{1_3})^T$ satisfies the relation  $(f_{1_1}^2-f_{1_2}^2)(f_{1_2}^2-f_{1_3}^2)(f_{1_3}^2-f_{1_1}^2)\neq 0$.
\end{example}
\begin{example}
    In the rigid body dynamics (\ref{par:RB}), note that the components of the $c$ vector cannot all be of the same sign for physical reasons. That's why in this case we have $ad_{f_1}^2f_0 \notin \mathcal{S}^1(f_0,f_1)(0)$ instead,  namely the Hermes-Sussmann necessary condition is not satisfied. As a result, the rigid body angular velocity dynamics  is not STLC from the origin with only a single pair of external torque applied. This is consistent with the results in \textnormal{\cite{kerai1995analysis}} echoed in \textnormal{\cite{coron2007control}}.
\end{example}

\section{Applications: The Sprott and Lorenz Systems}

Although a general condition for the single-input small-time local controllability
of systems  with $L \neq 0$ in the class $\mathcal{C}$ given in (\ref{class_S}) is still under investigation,  we're nevertheless able to derive specific necessary and sufficient conditions for the Sprott and Lorenz systems. It turns out that either one of these two systems is single-input STLC from the origin as long as the control vector avoids a specific subspace arrangement. The following two propositions are the main results of this section.

\paragraph{The Sprott system}
We have the following proposition about the single-input small-time local controllability of the Sprott system

\begin{proposition}\label{prop:Sprott-STLC=m1}
    The Sprott system (\ref{eq:SprottSystem-Delta-dynamics}) is small-time locally controllable from the origin with a single-input  control vector $f$ if and only if $f \notin Span\{\mathbf{1}\}\cup Span\{\mathbf{1}\}^{\perp}$
\end{proposition}

Before we proceed to the proof of this proposition, we make the following remark to highlight the conditions for single-input accessibility versus small-time local controllability from the origin of the Sprott system.
\begin{remark}
    In Proposition \ref{prop:propositionSubClass_Sprott} we've shown that a necessary and sufficient condition for the Sprott system to be single-input locally strongly  accessible from the origin is that the control vector $f$ not be along the diagonal $Span\{\mathbf{1}\}$. What Proposition \ref{prop:Sprott-STLC=m1} tells us is that if one wants to achieve small-time local controllability,  then the orthogonal complement of the diagonal namely $Span\{\mathbf{1}\}^{\perp}$ must be further avoided.
\end{remark}
\begin{proof}
    Note that $det\big\{f,Lf,L^2f \big\} = -\frac{1}{2}(f^T\mathbf{1})(f^THf)$, where $H$ is the hessian of the degenerate quadratic form $x^2+y^2+z^2-(yz+zx+xy)$ whose isotropic cone is the subspace $ker(H)=Span\{\mathbf{1}\}$. Thus if $f \notin Span\{\mathbf{1}\}\cup Span\{\mathbf{1}\}^{\perp}$ then the linearization at the origin is controllable, and the system is STLC. Conversely, on the one hand, if $f \in Span\{\mathbf{1}\}$ then the system is not even accessible per Proposition \ref{prop:propositionSubClass_Sprott}. On the other hand, if $f \in Span\{\mathbf{1}\}^{\perp}$ then a direct computation shows that the augmented linear system $\delimitershortfall=0pt
\setlength{\dashlinegap}{2pt}
\left[\begin{array}{ccc:c}
f & Lf &L^2f& \Phi(f)
\end{array}
\right]$ is inconsistent. As a result, we have $[[f_0,f],f] = 2\Phi(f) \notin \mathcal{S}^1(f_0,f)(0) = Span\big\{f,Lf,L^2f \big\}$, where the latter equality holds by Cayley-Hamilton. Therefore the system does not satisfy the Hermes-Sussmann necessary condition  and cannot be STLC from the origin.
\end{proof}

\paragraph{The Lorenz system}
We have the following proposition about the single-input small-time local controllability of the Lorenz system.

\begin{proposition}\label{prop:Lorenz-STLC=m1}
   Consider the Lorenz system given by (\ref{eq:Lorenz63-Delta-dynamics}), and define the constants $\mathfrak{s} = \beta^2-(\sigma+1)\beta + \sigma(1-\rho)$ and $ ~ \mathfrak{d}^2 = 4\rho\sigma+(\sigma-1)^2$. Furthermore, define the vectors $v_{\pm}=(1,\frac{\sigma-1 \pm \mathfrak{d}}{2\sigma},0)^T$ $w_{\pm}=(\frac{1-\sigma \pm \mathfrak{d}}{2\rho},1,0)^T$, and finally define the subspaces $V_{3+}= Span\{ \mathbf{e}_3,v_+\}, ~V_{3-}=Span\{ \mathbf{e}_3,v_-\}, ~W_{3+}=Span\{ \mathbf{e}_3,w_+\}$ and $W_{3-}=Span\{\mathbf{e}_3,w_-\}$. Under the assumption that $\mathfrak{s}\neq 0$, the Lorenz system is single-input small-time locally controllable from the origin with control vector $f$ if and only if   $f \notin   V_{3+}\cup V_{3-} \cup W_{3+} \cup W_{3-} \cup Span\{\mathbf{e}_3\}^\perp$.
\end{proposition}
The argument of the proof of this proposition is similar to that of Proposition \ref{prop:Sprott-STLC=m1} with the exception of a few details. Nevertheless, we present it below for the sake of completeness. 

\begin{proof}
    We note that $det\big\{f,Lf,L^2f \big\} = \frac{1}{2}\mathfrak{s}(f^T\mathbf{e}_3)(f^THf)$, where $H$ is the hessian of the degenerate quadratic form $\rho x^2-\sigma y^2 +(\sigma-1)xy$ whose isotropic cone is precisely the subspace arrangement $E=V_{3+}\cup V_{3-} \cup W_{3+} \cup W_{3-}$. Thus if $f\notin E\cup Span\{\mathbf{e}_3\}^\perp$ then the linearization at the origin is controllable and therefore the system is STLC. Conversely, if $f \in Span\{\mathbf{e}_3\}^\perp$, then $\mathcal{S}^1(f_0,f)(0) = Span\{ f,Lf,L^2f\} = Span\{\mathbf{e}_3\}^\perp$, and $\Phi(f) \in Span\{ \mathbf{e}_3\}$ since $\Phi(\mathbf{e}_1)=\mathbf{0}=\Phi(\mathbf{e}_2)$ and $\Psi(\mathbf{e}_1,\mathbf{e}_2)=\mathbf{e}_3$. If $f\in E$, then a  direct computation shows that the corresponding augmented linear system is always inconsistent. Consequently,  noting that $[[f_0,f],f] = 2\Phi(f)$, it follows that the system cannot be STLC from the origin by  the Hermes-Sussmann necessary condition.
\end{proof}

\section{Conclusions}
\label{sec:conclusions}
In  this work we have studied the local controllability of an important class of affine quadratic control systems with non necessarily homogeneous drift. The interest in this class is driven by the ubiquity in science and engineering of some of its distinguished members, namely the Sprott system, the Lorenz system and the rigid body. We've derived a necessary and sufficient condition for strong accessibility that is amenable to computer implementation, and it is   reminiscent of  the  Kalman rank condition for linear systems. Our condition recovers  that of Crouch in the case of the rigid body, and it's therefore more general. We have also established that our condition is necessary and sufficient to conclude small-time local controllability in the rank $1$ underactuation case under some relatively mild  assumptions. We are still investigating ways to bypass these assumptions, and  to cover higher order ranks of underactuation for STLC. Finally, we paid particular attention to the single-input small-time local controllability of the Sprott and Lorenz systems for which we're able to give explicit and easy-to-check necessary and sufficient conditions. One of the main implications of our central result is that it provides not only  an intuitive  and easily implementable condition to ascertain accessibility or lack thereof, but it also quantifies  the "degree of reachability"  of the system. The controllability results in this work will certainly pave the way to address in a systematic way the question of controlling  bifurcation and  chaos in the Sprott and Lorenz systems in particular, and by extension in this larger and rich class of systems.

\section*{Acknowledgments}
We would like to thank M. D. Kvalheim for his valuable comments and suggestions. The work was partially supported by NSF grant  DMS-2103026, and AFOSR grants FA 9550-22-1-0215 and FA 9550-23-1-0400.
\bibliographystyle{siamplain}
\bibliography{references}
\end{document}